\documentclass[12pt]{amsart}

\usepackage{ amssymb }
\usepackage{a4wide,xcolor}
\usepackage[alphabetic]{amsrefs}
\usepackage[T1,T2A]{fontenc}
\usepackage[utf8]{inputenc}
\usepackage{url,hyperref}
\usepackage{amsfonts}
\usepackage{mathtools}
\usepackage{needspace}
\usepackage[pdftex]{graphicx}
\usepackage{datetime}
\usepackage{epigraph}
\usepackage{verbatim}
\usepackage{mathtools}
\usepackage{enumerate}
\usepackage[american]{babel}
\numberwithin{equation}{section}

\newcommand{\Z}{\mathbb{Z}}
\newcommand{\T}{\mathbb{T}}
\newcommand{\N}{\mathbb{N}}
\newcommand{\R}{\mathbb{R}}

\newcommand{\Cm}{\mathbb{C}}

\newcommand{\eps}{\varepsilon}

\DeclareMathOperator{\supp}{supp}

\renewcommand{\phi}{\varphi}

\newtheorem{Thm}{Theorem}[section]
\newtheorem{theorem}[Thm]{Theorem}

\newtheorem{corollary}[Thm]{Corollary}

\newtheorem{conjecture}[Thm]{Conjecture}
\newtheorem{problem}[Thm]{Problem}

\begin{document}

\thanks{
Aleksei Kulikov was supported by the VILLUM Centre of Excellence for the Mathematics of Quantum Theory (QMATH) with Grant No.10059.
The work of Miquel Saucedo and  Sergey Tikhonov is supported by
 PID2023-150984NB-I00 funded by MICIU/AEI/10.13039/501100011033/ FEDER, EU, 
the CERCA Programme of the Generalitat de Catalunya and the Severo Ochoa, and Mar\'ia de Maeztu
Program for Centers and Units of Excellence in R\&d (CEX2020-001084-M).
Miquel Saucedo is supported by  the Spanish Ministry of Universities through the FPU contract FPU21/04230.
 Sergey Tikhonov is supported
by 2021 SGR 00087. }

\author{Aleksei Kulikov}
\address{Aleksei Kulikov, University of Copenhagen, Department of Mathematical Sciences,
Universitetsparken 5, 2100 Copenhagen, Denmark,
\newline {\tt lyosha.kulikov@mail.ru} 
}

\author{Miquel Saucedo}
\address{Miquel  Saucedo,  Centre de Recerca Matemàtica\\
Campus de Bellaterra, Edifici C
08193 Bellaterra (Barcelona), Spain}
\email{miquelsaucedo98@gmail.com }

\author{Sergey Tikhonov}
\address{Sergey Tikhonov, 
ICREA, Pg. Lluís Companys 23, 08010 Barcelona, Spain;
Centre de Recerca Matem\`{a}tica\\
Campus de Bellaterra, Edifici C
08193 Bellaterra (Barcelona), Spain,
 and Universitat Autònoma de Barcelona, Spain.}
\email{stikhonov@crm.cat}


  %\sloppy 

 % \title[A Kahane--Katznelson--de Leeuw theorem for sparse spectra]
  %{A Kahane--Katznelson--de Leeuw theorem for sparse spectra\\Fourier coefficients of continuous functions with sparse spectrum}
\title[Fourier coefficients of continuous functions with sparse spectrum]
{Fourier coefficients of continuous functions with sparse spectrum}

  %A lacunary version of the Kahane--Katznelson--de Leeuw theorem

  \begin{abstract}  
  Let $(r_k)$ be an increasing sequence and $(w_k)$ a positive sequence. We study the following question: is it true that for every sequence $(a_k)$ satisfying $\sum_{k=0}^\infty |a_k|^2 w_k^2 < \infty$ there exists a function $f\in C(\T)$ such that $\hat{f}(2^k) = a_k$ and $\hat{f}(n) = 0$ for $n\notin \cup_k [2^k-r_k,2^k+r_k]$? We show that this is possible if and only if $\sup_{k\in\N}\sum_{n=[\log_2 r_k]}^k w_k^{-2} < \infty$.
\end{abstract}
  \maketitle
  \markboth{Aleksei Kulikov, Miquel Saucedo, and 
Sergey Tikhonov}{Fourier coefficients of continuous functions with sparse spectrum}

  \section{Introduction and main results}
In this paper we work with functions on the unit circle
$\T$, which we view as the interval $[-\frac{1}{2},\frac{1}{2}]$ with endpoints identified. We equip it with the standard Lebesgue probability measure.  For an integrable  function  $f$ on $\T$, we define its Fourier coefficients by
$$\hat{f}(n) = \int_{-\frac{1}{2}}^{\frac{1}{2}} f(t) e^{-2\pi i n t}dt.$$
Parseval’s theorem states that the $L^2$-norm of a function is equal to the $\ell^2$-norm of its Fourier coefficients:
$$\int_{-\frac{1}{2}}^{\frac{1}{2}} |f(t)|^2dt = \sum_{n=-\infty}^\infty |\hat{f}(n)|^2.$$

Since all continuous functions 
 on 
$\T$ 
 are bounded, they belong to 
$L^2(\T)$, and hence their Fourier coefficients are square-summable. The celebrated Kahane–Katznelson–de Leeuw theorem \cite{DKK} states that this is the only possible restriction on the decay of the Fourier coefficients of a continuous function.
More precisely, 
\begin{theorem}\label{DKK theorem}
Let $(a_n)\in \ell^2(\Z)$. Then there exists a function $f\in C(\T)$ such that
\begin{equation}\label{DKK inequality}
|\hat{f}(n)|\ge |a_n|,
\end{equation}
for all $n\in\Z$.
\end{theorem}
Our study is motivated by the following question: can one replace the inequality in \eqref{DKK inequality} with an equality?
This turns out not to be the case.
If we assume that 
$a_n= 0$ whenever 
$n$ is not a power of 
$2$, then for any function 
$f$ such that 
$|\hat{f}(n)| = |a_n|$
 we necessarily have
$\|f\|_{C(\T)}\asymp \| a_n\|_{\ell^1(\Z)}$.
  It is therefore natural to ask whether we can have $|\hat{f}(n)| = |a_n|$ (or even $\hat{f}(n) = a_n$) for \textit{many} $n\in \Z$.
 %We focus on Fourier series that are close to lacunary and consider the following problem.
 In other words, we study Fourier series with a sparse spectrum and address the following problem.
 
\begin{problem}\label{DKK problem}
Let $(a_k)_{k\in\N_0}$ be a square-summable sequence, and let $(r_k)_{k\in \N_0}$ be a sequence of positive numbers. Does there exist a function $f\in C(\T)$ such that $$\hat{f}(2^k) = a_k$$ and
$$
\qquad
\hat{f}(n) = 0,\quad
n\notin \bigcup_{k\geq 0} [2^k-r_k, 2^k+r_k]?$$
\end{problem}
%Necessary and sufficient condition in terms of $a_k$ and $r_k$ is likely hopeless, so we will ask if it holds for all sequences $a_k$ such that
Following Theorem \ref{DKK theorem}, we ask whether Problem \ref{DKK problem} admits a solution for any sequence
$(a_k)_{k\in\N_0}$ satisfying
$$\sum_{k=0}^\infty |a_k|^2 w_k^2 < \infty,$$
where $(w_k)_{k\in\N_0}$ is a given weight sequence.

Let us also comment on the related property of massiveness. A set $A\subset \Z$ is called {\it massive} if for any sequence $(a_n)_{n\in A}\in \ell^2(A)$ there exists a function $f\in C(\T)$ such that $\hat{f}(n) = 0, n\notin A$ and $|\hat{f}(n)| \ge |a_n|, n\in A$. The Kahane--Katznelson--de Leeuw theorem says that $A = \Z$ is a massive set. Kislyakov showed that $\N_0$ (see \cite{KislyakovZ}), $\cup_{k=1}^\infty [2^{2k}, 2^{2k+1}]$ and  also $\cup_{k=1}^\infty [2^k(1-\varepsilon),2^k(1+\varepsilon)], 0<\varepsilon<1$ (see \cite{KislyakovE}) are massive sets. On the other hand, %it is well known that  
 the set $A=\{2^n\}_{n \geq 0}$ is not massive, since, as mentioned before, any function with $\hat{f}(n)=0$ if $n \not \in A$ satisfies $\|f\|_{C(\T)}\asymp \| \hat{f}(n)\|_{\ell^1(A)}$. 

 Our second problem is then to determine whether there exists a \textit{sparse} set $A$ which is massive. More precisely, we formulate
 \begin{problem}
 \label{problem massive}
 For which non-decreasing sequences $(r_n)_{n\geq 0}$ is the set $$A= \bigcup_{n \geq 0} [2^n-r_n,2^n+r_n]$$ massive?
 \end{problem}

Our main results are the following two theorems, which give necessary and sufficient conditions on $(r_k)_k$ and $(w_k)_k$ for an affirmative answer to Problem \ref{DKK problem}.
For convenience, we set 
{$w_k = \infty$} for 
$k < 0$, and by 
$[x]$ we denote the greatest integer not exceeding 
$x$.

\begin{theorem}\label{necessary}
Let $0 < r_k < 2^{k-10}$ be an increasing sequence  and $(w_k)_{k\in\N_0}$ be a sequence of positive numbers. If
$$\sup_{k\in\N_0}\sum_{n=[\log_2 r_k]}^k \frac1{w_n^{2}} = \infty,$$
then there exists a sequence $(a_k)_{k\in \N_0}$ such that 
$$\sum_{k=0}^\infty |a_k|^2 w_k^2 < \infty$$
but there does not exist a function $f\in C(\T)$ such that $|\hat{f}(2^k)|\ge |a_k|$ and
$$\hat{f}(n) = 0,\quad n\notin \bigcup_k [2^k-r_k, 2^k+r_k].$$
\end{theorem}
\begin{theorem}\label{sufficient}
Let $0 < r_k < 2^{k-10}$ be an increasing sequence  and $(w_k)_{k\in\N_0}$ be a sequence of positive numbers. If
$$\sup_{k\in\N}\sum_{n=[\log_2 r_k]}^k \frac1{w_n^{2}} < \infty,$$
then for all sequences $(a_k)_{k\in\N_0}$ such that
$$\sum_{k=0}^\infty |a_k|^2 w_k^2 < \infty$$
there exists a function $f\in C(\T)$ such that $\hat{f}(2^k)=a_k$ and
$$\hat{f}(n) = 0, \quad n\notin \bigcup_k [2^k-r_k, 2^k+r_k].$$
Moreover, we have 
$$\|f\|_{C(\T)} \le C\sqrt{\sum_{k=0}^\infty |a_k|^2 w_k^2}$$
for some absolute constant $C$ depending only on the sequences $w_k$ and $r_k$.
\end{theorem}
{
In particular, the answer to Problem~\ref{DKK problem} is affirmative}
for all $(a_k)_{k\in\N_0}$ satisfying\\
$\sum_{k=0}^\infty |a_k|^2 w_k^2 < \infty$ if and only if

$$\sup_{k\in\N_0}\sum_{n=[\log_2 r_k]}^k \frac1{w_n^{2}} <\infty.$$
 \iffalse
\begin{corollary}
    Let $(w_k)_{k\in\N_0}$ and  $(r_k)_{k\in \N_0}$ be sequences of positive numbers with $0<r_k<2^{k-10}$. Then for any $(a_k)_{k\in\N_0}$ satisfying
$\sum_{k=0}^\infty |a_k|^2 w_k^2 < \infty$ there exists  $f\in C(\T)$ such that $\hat{f}(2^k) = a_k$ and
$
\hat{f}(n) = 0,$ $
n\notin \bigcup_k [2^k-r_k, 2^k+r_k],$
if and only if
$$\sup_{k\in\N_0}\sum_{n=[\log_2 r_k]}^k \frac1{w_k^{2}} <\infty.$$
\end{corollary}
\fi
Combining Theorem \ref{necessary} with Kislyakov's result on the massiveness of $\cup_{k=1}^\infty [2^k(1-\varepsilon),2^k(1+\varepsilon)], 0<\varepsilon<1$, we provide a complete answer to Problem \ref{problem massive}:
\begin{corollary} Let $0 < r_n<2^{n-10}$ {be an increasing sequence}. Then the set
    $A= \bigcup_{n \geq 0} [2^n-r_n,2^n+r_n]$ is massive if and only if there exists $\varepsilon>0$ such that  $r_n > \eps 2^n.$
\end{corollary}

\iffalse Note that the condition that we consider in our work is closely related to massiveness -- while we demand that $\hat{f}(n) = 0, n\notin A$ where $A = \bigcup_k [2^k-r_k, 2^k+r_k]$ we require the inequality $|\hat{f}(n)|\ge |a_n|$ to only hold if $n = 2^k$.\textbf{???} 
\fi
Finally, we give a few comments about the techniques we use. Our argument to obtain Theorem \ref{sufficient} shares some steps with the aforementioned result of Kislyakov in \cite{KislyakovE}. However, when we arrive at the key interpolation inequality \eqref{interpolation} we are able to establish it through a weak-type estimate at $p = 1$ by using a smoothed multiplier. We note that in \cite{KislyakovE} such a weak-type estimate does not hold and the needed interpolation inequality was obtained directly. Lastly, because of the lacunarity we are able to turn inequality $|\hat{f}(2^k)|\ge |a_k|$ into an equality $\hat{f}(2^k) = a_k$.

%{\color{red}COMMENT ON NORM ESTIMATE FOR A FUNCTION}

The structure of the paper is as follows. In Section \ref{section2} we prove Theorem \ref{necessary}. In Section \ref{section3} we reduce Theorem \ref{sufficient}
to the interpolation inequality \eqref{interpolation}, which is proved in Section \ref{section4}.
 %. In Section \ref{section4} we prove the interpolation inequality \eqref{interpolation}. 
 %Finally, in Section \ref{section5} we  discuss what our methods give in the general setup of Problem \ref{DKK problem}.	
Finally, in Section \ref{section5}, we discuss the implications of our methods for the general setup of Problem \ref{DKK problem}.
 
\section{Proof of Theorem \ref{necessary}}\label{section2}
First of all, let us assume that $w_k^{-2}$ is unbounded. This means that there exists an increasing sequence of indices $k_n$ such that $w_{k_n}^{-2} > 2^{n}$. Let us consider the sequence $a_k$ with $a_{k_n} = 1$ and $a_{k} = 0$, $k\neq k_n$. Clearly $\sum |a_k|^2 w_k^2 < \infty$ but in this case there does not even exist $f\in L^2(\T)$ such that $|\hat{f}(2^k)|\ge |a_k|$ as $(a_k)_k\notin \ell^2(\N_0)$. From this point on, we assume that $w_{k}^{-2}$ is bounded from above.
 
Assume that there is a continuous function $f$ such that 
\begin{equation}
\label{eq:spectrumoff}|\hat{f}(2^k)|\ge |a_k|\; \;\mbox{ and } \;\;
    \hat{f}(n) = 0, \; n\notin \bigcup_k [2^k-r_k, 2^k+r_k].
\end{equation}
Our first goal is to show that for any such function
\begin{equation}
\label{eq:intermerdiategoal}\sum_{k=2M}^{2N} |a_k| \le 4\|f\|_{C(\T)},\end{equation} where $M<N$ satisfies $10r_{2M} < 4^N$.

To prove this, we use the simple fact that for any
$g \in L^1(\T)$, we have
$$\left|\int_{-\frac{1}{2}}^{\frac{1}{2}} f(x)\overline{g(x)}dx\right| \le \|f\|_{C(\T)}\|g\|_{L^1(\T)}.$$
Let us now construct  functions
 $g,h\in L^1(\T)$ such that $\|g\|_{L^1(\T)} = \|h\|_{L^1(\T)}=1$ and 
$$\left|\int_{-\frac{1}{2}}^{\frac{1}{2}} f(x)\overline{g(x)}dx\right| \geq 2 \sum_{k=M}^N |a_{2k}|, \quad \left|\int_{-\frac{1}{2}}^{\frac{1}{2}} f(x)\overline{h(x)}dx\right| \geq 2 \sum_{k=M}^N |a_{2k-1}|.
$$
 
We choose the function $g$ to be a Riesz product
$$g(x) = \prod_{k=M}^N \Big(1+\cos\,(2\pi 4^k (x+\alpha_k))\Big)$$
for some natural numbers $M < N$ and real numbers $\alpha_k$.

First, we show that \begin{equation}
\label{eq:normg1}
\|g\|_{L^1(\T)}=1.\end{equation}
To see this, observe that since $g(x) \ge 0, x\in [-\frac{1}{2}, \frac{1}{2}]$, we have $\|g\|_{L^1(\T)} = \int_{-\frac{1}{2}}^{\frac{1}{2}} g(x)dx = \hat{g}(0)$. To compute $\widehat{g}(0)$, we write $g(x)$ in terms of complex exponentials 
\begin{equation}\label{g product}
g(x) = \prod_{k=M}^N \left(1 + \frac{1}{2}e^{2\pi i 4^k (x+\alpha_k)}+\frac{1}{2}e^{-2\pi i 4^k (x+\alpha_k)}\right),
\end{equation}
and observe that the term corresponding to $\widehat{g}(0)$ 
arises from multiplying
 only the  $1$'s in each bracket. Indeed, suppose we could obtain a constant in some other way, and let $K$ be the largest index for which we did not take $1$ from the $K$‑th bracket.
In this bracket, we took either $\frac{1}{2} e^{2\pi i 4^K (x + \alpha_K)}$ or $\frac{1}{2} e^{-2\pi i 4^K (x + \alpha_K)}$. However, even if in all previous brackets we took powers with the opposite sign, the coefficient of $2\pi i x$ in the exponent is at least
$$4^K - \sum_{k=M}^{K-1}4^k > 4^k - \sum_{k=-\infty}^{K-1}4^k =  \frac{2}{3}4^K$$
in absolute value, in particular it can not be zero. 
Thus, the only way to obtain a constant term is to take $1$’s from all brackets.
Therefore, $\|g\|_{L^1(\T)} = 1$, that is, \eqref{eq:normg1} is proved.

Second, recalling that $4^M > 10r_{2N}$, we claim that, under this assumption, %we have
\begin{equation}\label{fg inner product}
\int_{-\frac{1}{2}}^{\frac{1}{2}} f(x)\overline{g(x)}dx = \frac{1}{2}\sum_{k=M}^N \hat{f}(4^k)e^{-2\pi i 4^k \alpha_k}
\end{equation} for any $f$ satisfying \eqref{eq:spectrumoff}.

By Parseval's theorem, we have
$$\int_{-\frac{1}{2}}^{\frac{1}{2}} f(x)\overline{g(x)}dx = \sum_{n\in \Z} \hat{f}(n)\overline{\hat{g}(n)}.$$

So, to prove \eqref{fg inner product}, it is enough to show that $$\hat{f}(n)\overline{\hat{g}(n)}=\begin{cases}
\frac{1}{2}\hat{f}(4^k)e^{-2\pi i 4^k \alpha_k}, \, &n = 4^k,\quad M \le k \le N,\\
0,&\mbox{ otherwise}.
\end{cases}$$

If $n\notin \bigcup_{k} [2^k-r_k, 2^k+r_k]$, then $\hat{f}(n)\overline{\hat{g}(n)}= 0$ as $\hat{f}(n) = 0$ for such $n$. So, it remains to show that
$${\hat{g}(n)}=\begin{cases}
\frac{1}{2}e^{2\pi i 4^k \alpha_k}, \, &n = 4^k,\;\; M \le k \le N\\
0,&\mbox{ otherwise}
\end{cases} \mbox{ \quad for }\quad  n \in \bigcup_{k} [2^k-r_k, 2^k+r_k].$$

To find $\hat{g}(n)$, we will again use the product formula \eqref{g product}. As in the computation of $\hat{g}(0)$, let us expand the product and let $K$ be the largest index such that we did not take $1$ from the $K$‑th bracket. If we took $\frac{1}{2} e^{2\pi i 4^K \alpha_K}$, then the coefficient $n$ of $2\pi i x$ in the exponent satisfies
$$4^K - \sum_{k=M}^{K-1} 4^k \le n \le 4^K + \sum_{k=M}^{K-1} 4^k,$$
while if we took $\frac{1}{2} e^{-2\pi i 4^K \alpha_K}$, it satisfies
$$-4^K - \sum_{k=M}^{K-1} 4^k \le n \le -4^K + \sum_{k=M}^{K-1} 4^k.$$
Note that in the second case $n$ is necessarily negative, so it will not be in our union of intervals. In the first case we have
$$\sum_{k=M}^{K-1}4^k < \sum_{k=-\infty}^{K-1} 4^k = \frac{4^K}{3}.$$
In particular, we have $\frac{2}{3}4^K < n < \frac{4}{3}4^K$. Hence, if $n \in [2^k - r_k, 2^k + r_k]$, then necessarily $k = 2K$, by our assumption that $r_k \le 2^{k-10}$.
If we took $1$’s from every other bracket, then we would obtain a contribution of $\frac{1}{2} e^{2\pi i 4^K \alpha_K}$ to the $4^K$‑th Fourier coefficient. Otherwise, let $S$ be the second largest index such that we did not take $1$ from the $S$‑th bracket.
For the index $n$ we therefore have
$$|n - 4^K| \ge 4^S - \sum_{k=M}^{S-1}4^k > \frac{2}{3}4^S \ge \frac{2}{3}4^M.$$
On the other hand, $r_{2K} \le r_{2N} < \frac{4^M}{10}$. Thus, $n$ is necessarily outside $[2^{2K} - r_{2K}, 2^{2K} + r_{2K}]$. So, the only non-zero Fourier coefficients of $g$ in $\bigcup_k [2^k-r_k, 2^k+r_k]$ are 
$$\hat{g}(4^k) = \frac{1}{2}e^{2\pi i 4^k \alpha_k}, \qquad M \le k \le N,$$ whence \eqref{eq:spectrumoff} follows.

Third, choose the numbers $\alpha_k\in \R$ so that $\hat{f}(4^k)e^{-2\pi i 4^k \alpha_k}=|\hat{f}(4^k)|$  for each $M \le k \le N$. In this case, we have
$$\int_{-\frac{1}{2}}^{\frac{1}{2}} f(x)\overline{g(x)}dx = \frac{1}{2}\sum_{k=M}^N |\hat{f}(4^k)| \ge \frac{1}{2}\sum_{k=M}^N |a_{2k}|.$$

Fourth, using a completely analogous argument with
$$h(x) = \prod_{k=M}^N \Big(1 + \cos(\pi 4^k(x+\beta_k))\Big),$$
we can show that 
$$\sum_{k=M}^N |a_{2k-1}| \le 2\|f\|_{C(\T)}.$$
Summing these, we obtain
$$\sum_{k=2M}^{2N} |a_k| \le 4\|f\|_{C(\T)},$$ and thus \eqref{eq:intermerdiategoal} is proved.

\vspace{2mm}
Finally, we construct a sequence $(a_k)$ such that
$\sum_{k=0}^\infty |a_k|^2 w_k^2<\infty$ but for which there exists no function $f$ satisfying \eqref{eq:spectrumoff}. To do so, let $M_s < N_s$ satisfy
\[
4^{M_s} > 10\, r_{2N_s} \quad \text{and} \quad \sum_{k=2M_s}^{2N_s} \frac1{w_k^{2}} > 16^s.
\]
To see that such $M_s$ and $N_s$ exist, note that, by the assumption that 
\(\sum_{j=[\log_2 r_k]}^k w_j^{-2}\) is unbounded, we can find $k_s$ such that
\[
\sum_{j=[\log_2 r_{k_s}]}^{k_s} \frac1{w_j^{2}} > 17^s.
\]
 Then, since we assumed that $w_k^{-2}$ is bounded,  we also have $$\sum_{j=10+[\log_2 r_{k_s}]}^{k_s+1} \frac1{w_j^{2}}>16^s, $$
 for $s$ large enough. Hence, we can choose $M_s$ and $N_s$ such that
\begin{align*}
    4. 5+[\log_2 r_{k_s}]/2&\leq M_s\leq 5+[\log_2 r_{k_s}]/2,\\
   k_s/2 &\leq N_s \leq k_s/2+0.5,
\end{align*}
% $$4. 5+[\log_2 r_{k_s}]/2\leq M_s\leq 5+[\log_2 r_{k_s}]/2 \mbox{ and } k_s/2\leq N_s \leq k_s/2+0.5,$$ 
respectively. {By increasing $k_s$, if necessary, we can also assume that $N_{s+1} > M_s$ so the intervals $[2M_s, 2N_s]$ for different $s$ are disjoint.} 
Consider the sequence defined by 
$$a_k = 2^{s}\frac{w_k^{-2}}{\sum_{n=2M_s}^{2N_s} w_n^{-2}},\quad k\in [2M_s, 2N_s]$$
and $a_k = 0$ otherwise. Then,
$$\sum_{k=0}^\infty |a_k|^2 w_k^2 = \sum_{s=0}^\infty 4^s \frac{1}{\sum_{n=2M_s}^{2N_s} w_n^{-2}} \le \sum_{s=0}^\infty 4^{-s} < \infty.$$\vspace{2mm} However, if there existed a function $f$ satisfying \eqref{eq:spectrumoff}, then from \eqref{eq:intermerdiategoal} we would deduce that
$%$\[
2^s = \sum_{k=2M_s}^{2N_s} |a_k| < 4 \|f\|_{C(\mathbb{T})}
$ %\]
 for all $s$, which is impossible. The proof is now complete.
%{\color{red}COMMENT ON BANACH--STEINHAUS AND SUBEXPONENTIAL CASE}
\section{Proof of Theorem \ref{sufficient}}
\label{section3}
Before starting the proof, we note that it is enough to consider the case where $r_k \to \infty$ and $(a_k)_k$ is finitely supported.

Indeed, we first note that the case of a uniformly bounded sequence $(r_k)$ is trivial. In this case,
$$\sup_{k \in \N} \sum_{n=[\log_2 r_k]}^k w_k^{-2}<\infty \implies (1/w_k)\in \ell ^2.$$ Therefore, the function
$$f(x) = \sum_{k=0}^\infty a_ke^{2\pi i 2^k x}$$
satisfies $\|f\|_{C(\T)}\le \|(a_k)\|_{\ell^1}\le \|(a_k w_k)\|_{\ell^2}\|1/w_k\|_{\ell^2}<\infty,$ whence the result follows. In what follows, we assume that $r_k \to \infty$.

Second, we show that it is enough to establish the result for finitely supported sequences $(a_k)_{k\in \N_0}$. Let $(a_k)_{k\in\N_0}$ be any sequence such that
$$\sum_{k=0}^\infty |a_k|^2 w_k^2 = t^2 < \infty$$ and assume that we have proved the theorem for finitely supported sequences.
Let $N_m$ be an increasing sequence of integers such that 
$$\sum_{k=N_m}^\infty |a_k|^2 w_k^2 \le t^24^{-m},\quad  m\ge 1.$$
For $m\ge 0$, consider sequences $a_{m, k} = a_k$ for $N_m \le k < N_{m+1}$ and $a_{m,k} = 0$ otherwise (where we put $N_0 = 0$). Each of these sequences is finitely supported, so we can find functions $f_m\in C(\T)$ such that $\hat{f}_m(2^k) = a_{m, k}$ and $\|f_m\|_{C(\T)} \le Ct2^{-m}$ (here, it is crucial that $C$ is independent of $N_m$). Thus, $f = \sum_{m=0}^\infty f_m$ is a continuous function as a uniformly convergent sum of continuous functions, $\hat{f}(2^k) = \sum_{m=0}^\infty a_{m,k} = a_k$ and $\|f\|_{C(\T)} \le 2Ct$, so we get the result by increasing $C$ by a factor of $2$. In the remainder of this section, we proceed in the spirit of Kahane--Katznelson--de Leeuw-type theorems; see, e.g., \cite{Nazarov,KislyakovE}.

\vspace{
2mm
}
\textbf{Step 0. Reformulation of the problem.}
We assume $a_k = 0$ for $k \ge N$. If all $a_k = 0$, then $f(x) = 0$ satisfies the conditions. Otherwise, we can assume without loss of generality that 
\begin{equation}
\label{eq:normalizationl2}
    \sum_{k=0}^\infty |a_k|^2 w_k^2 = 1
\end{equation}
by dividing each $a_k$ by $\sqrt{\sum_{k=0}^\infty |a_k|^2 w_k^2 }$. Consider the  finite-dimensional space $\mathcal{P}_{2^{N+1}}$ of functions on $\T$ whose frequencies lie between $1$ and $2^{N+1}$:
$$\mathcal{P}_{2^{N+1}} =\left\{v:\T\to \Cm\mid v(x) = \sum_{n=1}^{2^{N+1}} v_k e^{2\pi i nx}\right\}.$$
Define $g(x) = \sum_{k=0}^N a_k e^{2\pi i 2^k x}\in \mathcal{P}_{2^{N+1}}$ and consider the following two subsets of $\mathcal{P}_{2^{N+1}}$:
$$X = \left\{v\in \mathcal{P}_{2^{N+1}}\mid \hat{v}(n) = 0, n\notin \bigcup_{k\in\N_0} [2^k-r_k, 2^k+r_k],\,\, \lVert v\rVert_{C(\T)}\le C\right\}$$
and
$$Y = \{v\in \mathcal{P}_{2^{N+1}}\mid \hat{v}(2^k) = 0, k\in\N_0\},$$ where $C$ is a large enough constant.

Our ultimate goal is to show that $g\in X+Y$. Indeed, if $g = f+h, f\in X, h\in Y$, then $f$ satisfies all of the requirements.

For the sake of contradiction, we assume that $g \notin X + Y$.

\vspace{
2mm
}
\textbf{Step 1.}
Our first goal is to show that there exists a measure $\mu \in M(\T)$ such that
\begin{equation}
\label{eq:anotherinterm}
|\mu(g)|\geq 1,\quad\lVert \mu \rVert_{M(\T)} \leq \frac{1}{C} \quad\mbox{ and } \quad\hat{\mu}(n) = 0,\; n\in \bigcup_{k=0}^N [2^k-r_k, 2^k+r_k]\backslash \{2^k\},\end{equation} where we recall that for $\nu \in M(\T)$ we define its Fourier coefficients by
$$\hat{\nu}(n) = \overline{\nu(e^{2\pi i n t})} = \int_{-\frac{1}{2}}^{\frac{1}{2}} e^{-2\pi i n t} d\nu(t).$$
Note that both $X$ and $Y$ are absolutely convex, that is, they are convex and closed under multiplication by any complex number $z$ with $|z| = 1$.  Hence, $X+Y$ is also absolutely convex. If $g \notin X + Y$, then by the complex hyperplane separation theorem, there exists a linear functional $F \in (\mathcal{P}_{2^{N+1}})^*$ such that
\[
|F(v)| \le 1 \quad \text{for all } v \in X + Y, \quad \text{and} \quad |F(g)| \ge 1.
\] Note that we restricted ourselves to the finite-dimensional space $\mathcal{P}_{2^{N+1}}$ {only} for the purpose of applying the hyperplane separation theorem without any additional conditions to verify.

We use the functional $F$ in order to show that there exists a measure $\mu \in M(\T)$ which satisfies \eqref{eq:anotherinterm}.
We note that  $F(v) = 0$ for all $v\in Y$, equivalently,  $F(e^{2\pi i n x}) = 0$ for all $1 \le n \le 2^{N+1}$ 
that are not powers of $2$. Indeed, if $F(e^{2\pi i n x})\neq 0,$ then we could choose a constant
$A$ such that $|F(Ae^{2\pi i n x})| = |AF(e^{2\pi i n x})| > 1$. Since $Ae^{2\pi i n x} \in Y$ we arrive at a contradiction.

Next, we consider the following subspace of $\mathcal{P}_{2^{N+1}}:$ 
$$V= \left\{v\in \mathcal{P}_{2^{N+1}}\mid \hat{v}(n) = 0, n\notin \bigcup_{k\in\N_0} [2^k-r_k, 2^k+r_k]\right\},$$
equipped with the $C(\T)$-norm.

Since $X$ is the ball of radius $C$ in $V$, $F$ is a linear functional on $V$ with norm at most $\frac{1}{C}$. By the Hahn--Banach theorem, there exists a linear functional $\mu$ on $C(\mathbb{T})$ of norm at most $\frac{1}{C}$ which coincides with $F$ on $V$. Since $C(\mathbb{T})^* = M(\mathbb{T})$, the space of complex Borel measures on $\mathbb{T}$, $\mu$ can be identified with a measure $\overline{\mu} \in M(\mathbb{T})$ in the sense that
\[
\mu(f) = \int_{-\frac{1}{2}}^{\frac{1}{2}} f \, d\overline{\mu}.
\]
Using the fact that  $\mu$ coincides with $F$ on $V$, we derive that  
$$\hat{\mu}(n) = 0, \quad n\in \bigcup_{k=0}^N [2^k-r_k, 2^k+r_k]\backslash \{2^k\}$$ and $|\mu(g)|\geq 1$. This concludes the proof of \eqref{eq:anotherinterm} and Step 1.

\vspace{
2mm
}
\textbf{Step 2.} We now show that \eqref{eq:anotherinterm} leads to a contradiction by examining the properties of the Fourier multiplier operator \(M\) on \(M(\mathbb{T})\), defined by
\begin{equation}\label{abstract multiplier}
(M \nu)(t) = \sum_{n=1}^{\infty} m(n)\, \hat{\nu}(n)\, e^{2\pi i n t}.
\end{equation}
In order to define \(m(n)\), we introduce a sequence \((\lambda_s)_{s} \subset \mathbb{N}\) as follows: \(\lambda_0 = 0\) and
\[
r_{\lambda_{s+1}-1} < 2^{\lambda_s} \le r_{\lambda_{s+1}}.
\] Such a sequence exists, since by assumption \((r_k)\) is increasing and tends to infinity. Also, from  $r_k \le 2^{k-10}$ we deduce that $\lambda_{s+1} > \lambda_s$ and  $r_{\lambda_{s+1}} \ge 2^{10}r_{\lambda_s}$.  We define
$$m(n) = \begin{cases}\frac{\psi\left(\frac{n-2^k}{r_{\lambda_s}}\right)}{w_k},& 2^k-r_k \le n \le 2^k+r_k\,\,\,\, \text{for}\,\,\,\,\lambda_s \le k < \lambda_{s+1}\,\, \,\,\text{and}\,\,\,\,k\le N,\\
0, & \text{otherwise,}
\end{cases}$$
where $\psi$ is a fixed non-negative $C_0^\infty(\R)$ function supported on $[-1, 1]$ with $\psi(0) = 1$. Note that since $m(n)$ is non-zero for only finitely many $n$, the sum in \eqref{abstract multiplier} converges  for any measure $\nu$ and defines a smooth function.

In order to obtain a contradiction, we will use the facts that for some absolute constants $0<A,\delta<\infty$, we have
\begin{align}
  \label{eq:l2normMmu} &1\leq \|M\mu\|_{L^2(\T)};\\
     \label{eq:l43normMmu} &\delta \| M\mu\|_{L^2(\T)}\leq \|M\mu\|_{L^{4/3}(\T)};\\
    \label{interpolation}
&\|M\nu\|_{L^{4/3}(\T)} \le A \|\nu\|_{M(\T)}^{1/2} \|M\nu\|_{L^2(\T)}^{1/2}
 \quad{\text{for all}}\quad\nu \in M(\T).
\end{align}
Assuming for a moment that these three properties have been established, we obtain a contradiction. Indeed, by applying \eqref{interpolation} to $\nu=\mu$ and using  \eqref{eq:l43normMmu}, we obtain
$$\delta \| M\mu\|_{L^2(\T)}\le\|M\mu\|_{L^{4/3}(\T)} \le A \|\mu\|_{M(\T)}^{1/2} \|M\mu\|_{L^2(\T)}^{1/2} \leq \frac{A}{C^{\frac12}} \|M\mu\|_{L^2(\T)}^{1/2},$$ where the last inequality follows from \eqref{eq:anotherinterm}.  This is impossible if \(C > \frac{A^2}{\delta^2}\), in view of \eqref{eq:l2normMmu}. As a consequence, we deduce that $g\in X+Y$ for $C$ large enough and complete the proof.

All that remains is to establish  \eqref{eq:l2normMmu}--\eqref{interpolation}.

To prove \eqref{eq:l2normMmu}, consider the function $\tilde{g}(t) = \sum_{k=0}^N a_k w_k e^{2\pi i 2^k t}$, and observe that, by Parseval's theorem,
$$\mu(g) = \sum_{k=0}^N a_k\overline{\hat{\mu}(2^k)} = \sum_{k=0}^N (a_k w_k) \frac{\overline{\hat{\mu}(2^k)}}{w_k}=\int_{-\frac12}^{\frac12} \tilde{g}(t) \overline{M\mu(t)}dt,$$
where in the last step we used that $\psi(0) = 1$ so that $m(2^k) = \frac{1}{w_k}$ for $0 \le k \le N$.
Then, recalling that \(|\mu(g)| \geq 1\) and applying the Cauchy--Schwarz inequality, we obtain
\[
1 \leq \|M\mu\|_{L^2(\mathbb{T})} \, \|\tilde{g}\|_{L^2(\mathbb{T})},
\]
which, by the normalization condition \eqref{eq:normalizationl2} (i.e., \(\|\tilde{g}\|_{L^2(\mathbb{T})} = 1\)), yields \eqref{eq:l2normMmu}.

Next, to prove \eqref{eq:l43normMmu}, we recall that by definition of $m(n)$ and \eqref{eq:anotherinterm}, we have $$m(n) =0, n \not \in \bigcup_{k=0}^N [2^k-r_k, 2^k+r_k]\quad\mbox{ and } \quad \hat{\mu}(n) = 0, n\in \bigcup_{k=0}^N [2^k-r_k, 2^k+r_k]\backslash \{2^k\}.$$  As a consequence,  $M\mu$ is a lacunary series, that is,
$M\mu(t) = \sum_{k=0}^N b_k e^{2\pi i 2^k t}$ with  some $b_k\in\mathbb{C}$. For such series it is well known that all of its $L^p$-norms are equivalent for all $0 < p < \infty$. In particular, $\|M\mu\|_{L^{4/3}(\T)}\ge \delta \|M\mu\|_{L^2(\T)}$ for some absolute constant $\delta$. The proof of \eqref{eq:l43normMmu} is concluded.
 
 It remains only to establish \eqref{interpolation}. Since the proof is longer, it will be presented in the next section.
\section{Proof of the interpolation inequality \eqref{interpolation}}\label{section4}

\iffalse
I DONT THINK IT IS NECESSARY to writE this
First, we notice that it is enough to establish \eqref{interpolation} for measures represented by smooth functions. Indeed, for $K\in \N$ consider the Fejér kernel 
$$F_K(x) = \frac{1}{K} \left(\frac{1-\cos(2\pi K x)}{1-\cos(2\pi x)}\right)=\sum_{n=-K}^K \left(1-\frac{|n|}{K}\right)e^{2\pi i n x}$$
and for the measure $\nu$ let $\nu_K$ be the convolution of $\nu$ with $F_K$. Since $F_K\in L^1(\T)$ with $\|F_K\|_{L^1(\T)} = 1$, we have $\|\nu_K\|_{M(\T)} \le \|\nu\|_{M(\T)}$. On the other hand, since $\hat{\nu}_K(n) = \hat{\nu}(n)\hat{F}_K(n)$, $\nu_K$ is a smooth (in fact, even real analytic) function on $\T$. For each fixed $n$ we also have $\hat{\nu}_K(n) \to \hat{\nu}(n)$ since $\hat{F}_K(n)\to 1$. Thus, we have $M\nu_K\to M\nu$ uniformly (here, we used that $M$ has only finitely many non-zero terms and applied convergence to each of them separately). Therefore, $\|M\nu_K\|_{L^{4/3}(\T)}\to \|M\nu_K\|_{L^{4/3}(\T)}$ and $\|M\nu_K\|_{L^2(\T)} \to \|M\nu_K\|_{L^2(\T)}$. Passing to the limit $K\to \infty$ in 
$$\|M\nu_K\|_{L^{4/3}(\T)} \le A \|\nu_K\|_{M(\T)}^{1/2} \|M\nu_K\|_{L^2(\T)}^{1/2}\le A\|\nu\|^{1/2}_{M(\T)}\|M\nu_K\|_{L^2(\T)}^{1/2}$$
 we get \eqref{interpolation}.
\fi
In this section, we will make use of the Littlewood--Paley square function, which we now recall for completeness.
Let $\phi:\R\to\R$ be a smooth function supported on $[\frac{1}{4}+\frac{1}{100},	1-\frac{1}{100}]$ such that 
$$\sum_{k\in\Z} \phi(2^k x) = 1,\,\, x>0, \quad\mbox{ and } 
\quad
\phi(x)=1, \,\,x \in [1/2-1/{300},1/2+{1}/{300}].$$
Define the Fourier multiplier $P_k$ for $k\in \Z$ as follows:

$$P_k f(t)=\begin{cases}
      \sum_{n\in \Z} \phi(2^{-k} n) \hat{f}(n)e^{2\pi i n t}, &k>0,\\
      \sum_{n\in \Z} \phi(-2^{k} n) \hat{f}(n)e^{2\pi i n t}, &k<0,\\
      \hat{f}(0), &k=0.
\end{cases}$$
Note that, for any trigonometric polynomial $f$, we have $f(t) = \sum_{k\in \Z} P_k f(t)$ . We define the square function $Sf$ via
$$Sf(t) = \left(\sum_{k\in\Z} |P_k f(t)|^2\right)^{1/2}.$$
The Littlewood--Paley theory \cite[Theorem 6.1.2]{Grafakos} states that for \(1 < p < \infty\), the square function \(Sf\) has a norm comparable to that of \(f\), that is, there exist constants \(0 < c_p < C_p < \infty\) such that
\[
c_p \, \|f\|_{L^p(\mathbb{T})} \le \|Sf\|_{L^p(\mathbb{T})} \le C_p \, \|f\|_{L^p(\mathbb{T})}.
\]

We will apply this to $f = M\nu$ for $p = \frac{4}{3}$ and $p = 2$. So, to prove \eqref{interpolation}, it suffices to show that
$$\|SM\nu\|_{L^{4/3}(\T)} \le A'\|\nu\|^{1/2}_{L^1(\T)} \|SM\nu\|_{L^2(\T)}^{1/2}.$$
{\color{black}A direct computation with distribution functions shows that 
$$\|g\|_{L^{4/3}(\T)} \le C'\|g\|_{L^{1,\infty}(\T)}^{1/2}\|g\|_{L^{2,\infty}(\T)}^{1/2}\le C'\|g\|_{L^{1,\infty}(\T)}^{1/2}\|g\|_{L^{2}(\T)}^{1/2},$$
with some absolute constant $C'$, thus} \eqref{interpolation} follows from the weak-type estimate
\begin{equation}\label{weak}
\|SM\nu\|_{L^{1,\infty}(\T)} \le A'' \|\nu\|_{M(\T)}.
\end{equation}
\iffalse

We claim that for all $g\in L^2(\T)$ we have 

where $\|g\|_{L^{1,\infty}(\T)}= \sup_{x > 0} x |\{t: |g(t)|>x\}|$. Indeed, by Chebyshev inequality for $x > 0$ we have $|\{t: |g(t)| > x\}| \le \frac{\|g\|_{L^2(\T)}^2}{x^2}$ and by definition we have $|\{t: |g(t)| > x\}| \le \frac{\|g\|_{L^{1,\infty}(\T)}}{x}$. On the other hand, we have
$$\|g\|_{L^{4/3}(\T)}^{4/3} = \frac{4}{3}\int_0^\infty x^{1/3}|\{t: |g(t)| > x\}|dx.$$
Thus,
$$\|g\|_{L^{4/3}(\T)}^{4/3} \le \frac{4}{3}\int_0^\infty x^{1/3}\min\left(\frac{\|g\|_{L^{1,\infty}(\T)}}{x},  \frac{\|g\|_{L^2(\T)}^2}{x^2}\right)dx=6\|g\|_{L^{1,\infty}(\T)}^{2/3}\|g\|_{L^2(\T)}^{2/3}.$$
Raising this to the power $\frac{3}{4}$ we get the desired estimate.
\fi

As mentioned in the introduction, such a weak-type inequality is possible only because the multiplier \(M\) is constructed using a smooth cutoff function. In \cite{KislyakovE}, no such weak-type estimate is available, and the author instead proves the interpolation inequality \eqref{interpolation} directly. Moreover, the exponent on \(\|\nu\|_{M(\mathbb{T})}\) was not \(\tfrac{1}{2}\), but rather some \(0 < \alpha < \tfrac{1}{2}\).

Instead of the square function \(Sf\), we consider the closely related operator \(T\), which maps \(f\) to the vector
\[
Tf(t) = \bigl(P_k f(t)\bigr)_{k \in \mathbb{Z}},
\]
viewed as an element of \(\ell^2(\mathbb{Z})\). In other words, \(T\) maps a function \(f: \mathbb{T} \to \mathbb{C}\) to a function \(Tf: \mathbb{T} \to \ell^2(\mathbb{Z})\). We clearly have \(\|Tf(t)\|_{\ell^2(\mathbb{Z})} = |Sf(t)|\), but the advantage is that \(T\) is a linear operator. Hence, it suffices to establish the weak-type inequality for the vector-valued operator \(TM\).

Note that each coordinate of \(TM\) is a Fourier multiplier \(P_k M\), which corresponds, on the function side, to convolution with some function on \(\mathbb{T}\). Thus, \(TM\) as a whole is a vector-valued convolution operator with kernel \(K : \mathbb{T} \to \ell^2(\mathbb{Z})\). In particular, classical sufficient conditions guaranteeing weak-type bounds for such operators (see, e.g., \cite[Theorem 5.6.1]{Grafakos}) can be applied. There are four such conditions, depending on a parameter \(D \ge 0\).
\begin{enumerate}
\item There exists some $K_0\in \ell^2(\Z)$ such that $$\lim_{\eps\to 0} \left\|\int_{\eps < |y| < 1/2} K(y)dy - K_0\right\|_{\ell^2(\Z)} = 0.$$
\item The convolution operator with kernel \(K\) extends to a continuous map from \(L^2(\mathbb{T})\) to \(L^2(\mathbb{T}, \ell^2(\mathbb{Z}))\) with norm at most \(D\).
\item For all $-\frac{1}{2} < y < \frac{1}{2}$, we have $\|K(y)\|_{\ell^2(\Z)} \le \frac{D}{|y|}$. %\textbf{ I think that here all that matters is that this is finite}
\item For all $-\frac{1}{4} < y < \frac{1}{4},$ we have $$\int_{2|y| < |x| < \frac{1}{2}} \|K(x-y)-K(x)\|_{\ell^2(\Z)}dx \le D.$$
\end{enumerate}

If these four conditions are satisfied, then the convolution operator with kernel \(K\) satisfies a weak-type inequality from \(L^1(\mathbb{T})\) to \(L^{1,\infty}(\mathbb{T}, \ell^2(\mathbb{Z}))\) with norm at most \(cD\), for some absolute constant \(c\) (depending only on the underlying space \(\mathbb{T}\)). Thus, it remains to verify these conditions with a uniform bound on \(D\) for \(TM\).

Note that the first condition is the  one that does not involve \(D\),  it is required only  to ensure that the operator is well-defined. In our situation, since the multiplier \(M\) is finitely supported, \(K\) is represented by a smooth function with values in a finite-dimensional subspace of \(\ell^2(\mathbb{Z})\), in particular, the integral converges absolutely.

The second condition follows from the facts that $M$ is bounded on $L^2(\mathbb{T})$ and $T$ is bounded from $L^2(\T)$ to $L^2(\T,\ell^2(\Z))$. Indeed, the boundedness of $M$ follows from the fact that its symbol $m(n)$ is uniformly bounded by $\sup_k w_k^{-1}\le C^{\frac12}$; and the boundedness of $T$ follows from Parseval's theorem.
\iffalse use Parseval's theorem. To show that a Fourier multiplier is bounded from $L^2(\T)$ it is necessary and sufficient to show its symbol is uniformly bounded. Since $m(n)$, the symbol of $M$ is boundedFor the operator $TM$ each coordinate is multiplied by some value of $\psi$,  some value of $\phi$ and $\frac{1}{w_k}$ for some $k$. Since $\psi$ and $\phi$ are smooth and compactly supported, first two coefficients are uniformly bounded. By applying the assumption $\sum_{[\log_2 r_k]}^k w_n^{-2} \le C$ with $r_k\le 2^{k-10}$ we can leave only $n=k$ term and see that $w_k^{-2}\le C$, hence $\frac{1}{w_k} \le \sqrt{C}$ is also uniformly bounded. Thus, our operator is uniformly bounded from $L^2(\T)$ to $L^2(\T, \ell^2(\Z))$.\fi

For the remaining two conditions, we observe that $K=(K_k)_{k \in \Z}$ with
$$K_k(x) = \sum_{n \in \Z} \hat{P}_k(n) m(n)=\frac{1}{w_{k-1}}\sum_{n\in \Z}\phi(2^{-k}n)\psi\left(\frac{n-2^{k-1}}{r_{\lambda_s}}\right)e^{2\pi i n x},$$ where  $\lambda_s < k \le \lambda_{s+1}$ and $1 \le k \le N+1$; and $K_k=0$ otherwise.
Note that we used here that \(r_{\lambda_s} \le r_{k-1} \le 2^{k-11}\), so that the only part of the multiplier \(M\) appearing in \(K_k\) after applying \(P_k\) is 
\[
\frac{\psi\bigl(\frac{n - 2^{k-1}}{r_{\lambda_s}}\bigr)}{w_{k-1}}.
\]
Additionally, since any \(n\) for which 
\(\psi\bigl(\frac{n-2^{k-1}}{r_{\lambda_s}}\bigr) \neq 0\) lies in the interval \(2^{k-1}-r_{\lambda_s} < n < 2^{k-1}+r_{\lambda_s}\), and recalling that \(\phi(2^{-k} n) = 1\) for such \(n\), we can simplify further as follows:
\[
K_k(x) = \frac{1}{w_{k-1}} \sum_{n \in \mathbb{Z}} \psi\Bigl(\frac{n-2^{k-1}}{r_{\lambda_s}}\Bigr) e^{2\pi i n x}
= \frac{1}{w_{k-1}} e^{2\pi i 2^{k-1} x} \sum_{n \in \mathbb{Z}} \psi\Bigl(\frac{n}{r_{\lambda_s}}\Bigr) e^{2\pi i n x}.
\]

Let us denote this last sum by $R_{s}(x)$. We claim that for all $s\ge 0$ and all $0 < y < \frac{1}{2}$ we have
\begin{equation}\label{R bound}
|R_{s}(y)|\le \frac{C r_{\lambda_s}}{(1+(r_{\lambda_s}|y|)^2)}
\end{equation}
and for all $-\frac{1}{2} < y < \frac{1}{2}$ and for all $x$ with $2|y| < |x| < \frac{1}{2}$ we have
\begin{equation}\label{R-R bound}
|R_s(x-y)-R_s(x)|\le \frac{Cr_{\lambda_s} \min(r_{\lambda_s} |y|, 1)}{(1+(r_{\lambda_s}|x|)^2)},
\end{equation}
where $C$ is an absolute constant (depending only on the choice of the smooth function $\psi$).

For \eqref{R bound}, if $|y| \le \frac{1}{r_{\lambda_s}}$ then we simply bound $|e^{2\pi i n x}| \le 1$ and derive that $|R_s(x)| \le 2 \sup_{t\in\R} \psi(t) r_{\lambda_s}$ since $\supp \psi \subset [-1, 1]$. If $|y| > \frac{1}{r_{\lambda_s}}$ then we will use summation by parts twice to get
$$R_s(y) = \sum_{n\in \Z} \left(\psi\left(\frac{n}{r_{\lambda_s}}\right)-2\psi\left(\frac{n-1}{r_{\lambda_s}}\right)+\psi\left(\frac{n-2}{r_{\lambda_s}}\right)\right)\frac{e^{2\pi i n y}-1}{(e^{2\pi i y}-1)^2}.$$
The first term here is \(\frac{1}{r_{\lambda_s}^2}\psi''(t)\) for some \(t \in \mathbb{R}\), which in particular does not exceed \(\frac{c}{r_{\lambda_s}^2}\). The second term is bounded by \(\frac{1}{|y|^2}\), from which we obtain the desired estimate, since there are at most \(4 r_{\lambda_s}\) non-zero terms in the sum.

For \eqref{R-R bound}, if $r_{\lambda_s} |y| \ge 1$ then we simply use bound \eqref{R bound} for both $R_s(x-y)$ and $R_s(x)$ and we have to multiply $C$ by at most $8$. Otherwise, we write
$$R_s'(x) = \sum_{n\in\Z} \psi\left(\frac{n}{r_{\lambda_s}}\right)2\pi i n e^{2\pi i n x} = 2\pi i r_{\lambda_s} \sum_{n\in\Z}  \frac{n}{r_{\lambda_s}}\psi\left(\frac{n}{r_{\lambda_s}}\right) e^{2\pi i n x},$$
so, repeating the previous argument with the smooth compactly supported function $t\psi(t)$, we obtain  $$|R_s'(x)| \le r_{\lambda_s}\frac{C r_{\lambda_s}}{(1+(r_{\lambda_s} |x|^2))}.$$ Since $R_s(x-y)-R_s(x) = yR_s'(t)$ for some $t$ between $x-y$ and $x$ we arrive at the desired estimate.

Next, we apply these bounds to estimate $K_k(x)$. From \eqref{R bound} we obtain
$$|K_k(y)| \le \frac{1}{w_{k-1}} \frac{C r_{\lambda_s}}{(1+(r_{\lambda_s} |y|)^2)}.$$
For $K_k(x-y) - K_k(x)$ we have
\begin{align*}
K_k(x-y) - K_k(x) = \frac{1}{w_{k-1}}e^{2 \pi i 2^{k-1}(x-y)}(R_{s}(x-y)-R_{s}(x)) \\+ \frac{1}{w_{k-1}}\left(e^{2\pi i 2^{k-1}(x-y)}-e^{2\pi i 2^{k-1} x}\right)R_{s}(x).
\end{align*}
Plugging in our bounds and using $\left|e^{2\pi i 2^{k-1}(x-y)}-e^{2\pi i 2^{k-1} x}\right| \le \min (4\pi, 2\pi 2^k |y|)$ and $r_{\lambda_s} \le r_k \le 2^k \le 2^{\lambda_{s+1}}$ we obtain
$$|K_k(x-y)-K_k(x)| \le \frac{1}{w_{k-1}} \frac{Cr_{\lambda_s} \min(1, 2^{k-1} |y|)}{(1+(r_{\lambda_s} |x|)^2)} \le \frac{C}{w_{k-1}} \frac{r_{\lambda_s}\min(1, 2^{\lambda_{s+1}}|y|)}{(1+(r_{\lambda_s} |x|)^2)}.$$

Note that we established these estimates for $1 \le k \le N+1$, for other $k$ we have $K_k(y) = 0$, so 
they hold trivially. 

When we combine all values of \(k\) to compute \(\|K(y)\|_{\ell^2(\mathbb{Z})}\) and \(\|K(x-y)-K(x)\|_{\ell^2(\mathbb{Z})}\), we first square the terms, sum them, and then take the square root. We will repeatedly use the simple inequality 
$
\sqrt{u+v} \le \sqrt{u} + \sqrt{v},
$ 
and split the sums over \(\lambda_s < k \le \lambda_{s+1}\) to obtain
$$\|K(y)\|_{\ell^2(\Z)} \le \sum_{s=0}^\infty \sqrt{\sum_{k=\lambda_s+1}^{\lambda_{s+1}}\frac{C^2}{w_{k-1}^2} \frac{r_{\lambda_s}^2}{(1+(r_{\lambda_s}|y|)^2)^2}}.$$

By our assumption and the definition of the sequence $\lambda_s$, we  see that the sums $\sum_{k=\lambda_s+1}^{\lambda_{s+1}} {1}/{w_{k-1}^2}$ are uniformly bounded in $s$. Thus,
$$\|K(y)\|_{\ell^2(\Z)} \le C'\sum_{s=0}^\infty \frac{r_{\lambda_s}}{1+(r_{\lambda_s}|y|)^2)}.$$
Note that this is the key reason why we required the intermediate step using the square function, without it we would instead have the sum $\sum_{k=\lambda_s+1}^{\lambda_{s+1}} \frac{1}{w_{k-1}}$ which is not necessarily bounded.

Let $s_y$ be such that $r_{\lambda_{s_y}} < \frac{1}{|y|} \le r_{\lambda_{s_y+1}}$ and split this sum accordingly. We have
$$\|K(y)\|_{\ell^2(\Z)} \le C'\sum_{s=0}^{s_y} r_{\lambda_s} + C'\sum_{s=s_y+1}^\infty \frac{1}{r_{\lambda_s}|y|^2}.$$
Since $r_{\lambda_{s+1}} \ge 2r_{\lambda_s}$ we can bound the first sum by at most twice its last term and the second sum by at most twice its first term to get
$$\|K(y)\|_{\ell^2(\Z)} \le 2C'\left(r_{\lambda_{s_y}}+\frac{1}{r_{\lambda_{s_y+1}}|y|^2}\right)\le \frac{4C'}{|y|},$$
as required.

Finally, applying the same argument to verify the fourth property, we obtain
$$\int_{2|y|}^{1/2} \|K(x-y)-K(x)\|_{\ell^2(\Z)} dx\le C'\int_{2|y|}^{1/2}\sum_{s=0}^\infty \frac{r_{\lambda_s}\min(1, 2^{\lambda_{s+1}}|y|)}{(1+(r_{\lambda_s}|x|)^2)}dx.$$
We will choose the same \(s_y\) such that 
$r_{\lambda_{s_y}} < \frac{1}{|y|} \le r_{\lambda_{s_y+1}}$
and split the sum according to \(s < s_y-1\), \(s = s_y-1\), \(s = s_y\), and \(s > s_y\).

For $s > s_y$ we have
$$\int_{2|y|}^{1/2}\sum_{s=s_y + 1}^\infty \frac{r_{\lambda_s}\min(1, 2^{\lambda_{s+1}}|y|)}{(1+(r_{\lambda_s}|x|)^2)}dx \le \int_{2|y|}^\infty \sum_{s=s_y+1}^\infty\frac{1}{r_{\lambda_s}|x|^2}dx = \frac{1}{2|y|}\sum_{s=s_y+1}^\infty \frac{1}{r_{\lambda_s}}.$$
The last sum is at most twice its first term, so we bound it by  $\frac{1}{|y|r_{\lambda_{s_y}+1}} \le 1$. For $s = s_y$ and   $s = s_y-1$
we have
\begin{align*}
    &\int_{2|y|}^{1/2} \frac{r_{\lambda_{s_y}}\min(1, 2^{\lambda_{s_y+1}}|y|)}{(1+(r_{\lambda_{s_y}}|x|)^2)}dx \le \int_0^\infty \frac{r_{\lambda_{s_y}}}{(1+(r_{\lambda_{s_y}}|x|)^2)}dx = \frac{\pi}{2},
    \\
    &\int_{2|y|}^{1/2} \frac{r_{\lambda_{s_y-1}}\min(1, 2^{\lambda_{s_{y}}}|y|)}{(1+(r_{\lambda_{s_{y}-1}}|x|)^2)}dx \le \int_0^\infty \frac{r_{\lambda_{s_{y}-1}}}{(1+(r_{\lambda_{s_{y}-1}}|x|)^2)}dx = \frac{\pi}{2}.
\end{align*}
\iffalse
$$\int_{2|y|}^{1/2} \frac{r_{\lambda_{s_y}}\min(1, 2^{\lambda_{s_y+1}}|y|)}{(1+(r_{\lambda_{s_y}}|x|)^2)}dx \le \int_0^\infty \frac{r_{\lambda_{s_y}}}{(1+(r_{\lambda_{s_y}}|x|)^2)}dx = \frac{\pi}{2}.$$
For $s = s_y-1$ we similarly have
$$\int_{2|y|}^{1/2} \frac{r_{\lambda_{s_y-1}}\min(1, 2^{\lambda_{s_{y}}}|y|)}{(1+(r_{\lambda_{s_{y}-1}}|x|)^2)}dx \le \int_0^\infty \frac{r_{\lambda_{s_{y}-1}}}{(1+(r_{\lambda_{s_{y}-1}}|x|)^2)}dx = \frac{\pi}{2}.$$\fi
Finally, for $s < s_y-1$ we estimate
\begin{align*}
    \int_{2|y|}^{1/2}\sum_{s=0}^{s_y-2} \frac{r_{\lambda_s}\min(1, 2^{\lambda_{s+1}}|y|)}{(1+(r_{\lambda_s}|x|)^2)}dx &\le \int_{0}^{\infty}\sum_{s=0}^{s_y-2} 2^{\lambda_{s+1}}|y|\frac{r_{\lambda_s}}{(1+(r_{\lambda_s}|x|)^2)}dx\\&=\frac{\pi |y|}{2}\sum_{s=0}^{s_y-2}2^{\lambda_{s+1}}.
    \end{align*}
Since $\lambda_s$ is increasing, the sum is at most twice its last term and we bound it by $$\pi |y| 2^{\lambda_{s_y-1}} \le \pi |y| r_{\lambda_{s_y}} \le \pi.$$ Overall, the integral is at most \(C'(2\pi + 1)\), that is, it is uniformly bounded, as required. This concludes the proof of \eqref{weak}, and consequently of \eqref{interpolation}.

\vspace{2mm}
\section{Removing the weight}
\label{section5}

In this section, we briefly discuss what our arguments yield for the general Problem~\ref{DKK problem}. If the sequence \((r_k)\) is uniformly bounded, then for every function \(f \in C(\mathbb{T})\) such that 
\(\hat{f}(n) = 0\) for \(n \notin \bigcup_k [2^k - r_k, 2^k + r_k]\), 
we must have 
$
\sum_k |\hat{f}(2^k)| < \infty.
$ 
Conversely, for every sequence \((a_k)_{k \in \mathbb{N}_0} \in \ell^1(\mathbb{N}_0)\), 
trivially,
 there exists  \(f \in C(\mathbb{T})\) such that 
\(\hat{f}(n) = 0\) for \(n \notin \bigcup_k [2^k - r_k, 2^k + r_k]\) and \(\hat{f}(2^k) = a_k\), that is, $f(x) = \sum_{k} a_ke^{2\pi i 2^k x}$.

Henceforth, we will assume that \(r_k \to \infty\). Define the sequence \((\lambda_s)_s\) as before, with 
\(\lambda_0 = 0\) and 
$r_{\lambda_{s+1}-1} < 2^{\lambda_s} \le r_{\lambda_{s+1}}$.

Given a sequence \((a_k)_{k \in \mathbb{N}_0}\) of complex numbers, it is natural to ask  whether there exists a sequence \((w_k)_{k \in \mathbb{N}_0}\) for which Theorem~\ref{sufficient} can be applied. It is not hard to see that this is possible if and only if 
\[
\sum_{s=0}^\infty \left(\sum_{k=\lambda_s}^{\lambda_{s+1}-1} |a_k|\right)^2 < \infty.
\] 
Hence, we obtain the following result.

\begin{theorem} 
\label{th:halfconj} Let $0 < r_k < 2^{k-10}$ be an increasing sequence tending to infinity. Define the sequence $\lambda_s$ by $\lambda_0 = 0$ and $r_{\lambda_{s+1}-1} < 2^{\lambda_s} \le r_{\lambda_{s+1}}$. Let $(a_k)_{k\in\N_0}$ be a sequence of complex numbers satisfying %. If
$$\sum_{s=0}^\infty \left(\sum_{k=\lambda_s}^{\lambda_{s+1}-1} |a_k|\right)^2 < \infty.$$
Then there exists a function $f\in C(\T)$ such that $\hat{f}(n) = 0, n\notin\bigcup_k [2^k-r_k, 2^k+r_k]$ and $\hat{f}(2^k) = a_k$.  
\end{theorem}
In fact, by carefully following the proof, one can even guarantee that 
\[
\|f\|_{C(\mathbb{T})}^2 \le C \sum_{s=0}^\infty \Biggl( \sum_{k=\lambda_s}^{\lambda_{s+1}-1} |a_k| \Biggr)^2,
\] 
for some absolute constant \(C\) which does not depend on the sequence \((r_k)\) (that is, it depends only on the underlying space \(\mathbb{T}\)).

\iffalse
Thus, $\ell^2$ condition in $s$ is sufficient  and $\ell^\infty$ condition in $s$ is necessary. If $r_k \ge \eps 2^k$ then $\ell^2$ condition is equivalent to $(\hat{f}(2^k))_{k\in\N_0}$ being in $\ell^2(\N_0)$ which is of course always necessary as $C(\T)\subset L^2(\T)$.
\fi
We conjecture that the condition in Theorem \ref{th:halfconj} is in fact necessary. More precisely:
\begin{conjecture}\label{conj}
Let $0 < r_k < 2^{k-10}$ be an increasing sequence, tending to infinity.  Let the sequence \((\lambda_s)\) be defined by $\lambda_0 = 0$ and $r_{\lambda_{s+1}-1} < 2^{\lambda_s} \le r_{\lambda_{s+1}}$. There exists a constant $C$ such that for all $f\in C(\T)$ such that $\hat{f}(n) = 0$ if $n\notin \bigcup_k [2^k-r_k, 2^k+r_k]$ we have
$$\sum_{s=0}^\infty \left(\sum_{k=\lambda_s}^{\lambda_{s+1}-1} |\hat{f}(2^k)|\right)^2 \le C\|f\|_{C(\T)}^2.$$
\end{conjecture}
It is worth noting that the argument in Section~2 shows that for any function \(f \in C(\mathbb{T})\) with 
\(\hat{f}(n) = 0\) for \(n \notin \bigcup_k [2^k - r_k, 2^k + r_k]\), we have

$$\sup_s \sum_{k=\lambda_s}^{\lambda_{s+1}-1} |\hat{f}(2^k)|\le 9\|f\|_{C(\T)}.$$
Note that if Conjecture
\ref{conj}
were true, then we would 
%be able to completely eliminate the weight $(w_k)_{k\in\N_0}$ from our results and
deduce that, for a given sequence $(a_k)_{k\in\N_0}$, there exists $f\in C(\T)$ with $\hat{f}(n) = 0, n\notin\bigcup_k [2^k-r_k, 2^k+r_k]$ and $\hat{f}(2^k) = a_k$ if and only if 
$$\sum_{s=0}^\infty \left(\sum_{k=\lambda_s}^{\lambda_{s+1}-1} |\hat{f}(2^k)|\right)^2 <\infty.$$

\vspace{5mm}

\subsection*{Acknowledgments} 
Part of this project was carried out while Aleksei Kulikov was participating in the Intensive Research Programme on Modern Trends in Fourier Analysis at Recerca de Matemàtica in Barcelona during May–June 2025. He would like to thank the host and organizers for their hospitlity. 

{The authors would like to thank Sergei Kislyakov for helpful discussions  regarding the results of the paper \cite{KislyakovE}}.

\end{document}